\documentclass[11pt,a4paper,notitlepage,final,onecolumn,oneside]{article}
\usepackage{exscale}
\usepackage{latexsym}
\usepackage{calc}
\usepackage{amssymb}
\usepackage{url}

\newcommand{\proof}{{\bfseries Proof:\quad}}

\begin{document}

\newtheorem{lemma}{Lemma}
\newtheorem{proposition}[lemma]{Proposition}
\newtheorem{theorem}[lemma]{Theorem}
\newtheorem{definition}[lemma]{Definition}
\newtheorem{hypothesis}[lemma]{Hypothesis}
\newtheorem{conjecture}[lemma]{Conjecture}
\newtheorem{problem}[lemma]{Problem}
\newtheorem{remark}[lemma]{Remark}
\newtheorem{example}[lemma]{Example}
\newtheorem{property}[lemma]{Property}
\newtheorem{corollary}[lemma]{Corollary}
\newtheorem{algorithm}[lemma]{Algorithm}


\title{On the nature of finite groups}
\author{Aleksandr Golubchik\\
        \small Bramscher Str. 57, 49088 Osnabrueck, Germany,\\
        \small e-mail: agolubchik@gmx.de}

\date{13.08.2004}
\maketitle

\begin{abstract}
The reality of the difficulties in an investigation of finite groups are
considered. It is shown that the consideration of symmetry properties of
the $k$-orbits that are obtained with an action of a finite group
$F=(V,\cdot)$ on Cartesian power $V^k$ gives a new view on the nature of
groups and simplifies some difficult properties of groups.

Using this representation it is obtained a simple proof of the W. Feit,
J.G. Thompson theorem: Solvability of groups of odd order.\bigskip


\end{abstract}


\section{Introduction}

The group theory was born as the permutation group theory and later was
abstracted to a group algebra with corresponding properties of a group
operation.

Let $F$ be an abstract group on a set $F$, then the group algebra is
equivalent to the action of $F$ on $F$. But this algebra generates also the
action of $F$ on $F^k$. It is clear that that action can have its own
properties which belong of course also to group properties, but those
properties lie out of group algebra.

Namely a consideration of such properties joined with a term
\emph{$k$-orbit theory} was first considered by author in
\cite{Golubchik}. A full text of $k$-orbit theory with applications is
planned to be published later. Here we consider how this theory leads to
a simple proof of the W. Feit, J.G. Thompson theorem: Solvability of
groups of odd order (it is known that the original text covers 255 pages
\cite{Gorenstein}).

Below under primitive permutation group we understand non-Abelian
primitive group.


\section{$n$-Orbit representation of finite groups}

Let $F$ be a finite group, $A<F$, $|F|/|A|=n$ and $\overrightarrow{L_n}$,
$\overrightarrow{R_n}$ be ordered sets of left, right cosets of $A$
in $F$. It is known that every transitive permutation representation of
$F$ is equivalent to a representation of $F$ given by \emph{$n$-orbits}
$X_n'=\{f\overrightarrow{L_n}:\, f\in F\}$ or
$X_n''=\{\overrightarrow{R_n}f:\, f\in F\}$. It is also known that $F$ is
homomorphic to its image $Aut(X_n')=Aut(X_n'')$ with the kernel of the
homomorphism equal to a maximal normal subgroup of $F$ contained in $A$.
Further we assume that a finite group is isomorphic to its permutation (=
$n$-orbit) representation.

A maximal by inclusion subgroup $A$ of a finite group $F$ that contains no
normal subgroup of $F$ we call a \emph{md-stabilizer} of $F$ and a
corresponding representation of $F$ we call a \emph{md-representation} or
a \emph{md-group}. A md-stabilizer $A$ of a finite group $F$ defines a
minimal degree permutation representation of $F$ in a family of
permutation representations of $F$ produced with subgroups of $A$. A
non-minimal degree representation of $F$ we call a
\emph{nmd-representation} or a \emph{nmd-group}. A finite group can have
many (non-conjugate) md-stabilizers.

A special interest is presented by a permutation representation of the
lowest degree or a \emph{ld-representation} or a \emph{ld-group},
because a $n$-orbit of such representation (of degree $n$) and its
$k$-projections ($k$-orbits) contain all specific symmetry properties
which describe a finite group. This representation not always can be
represented with permutations of cosets of some subgroup of finite group.
As distinct from ld-representation, $k$-orbits of other permutation
representations can have additional symmetry properties that are not
specific for a finite group, but describe properties of corresponding
permutation group as for example the automorphism group of Petersen graph.
Some properties of ld-representation (and also other representations) are
on principle combinatorial and cannot be interpreted with group algebra.
For example, a ld-representation of a direct product and only that
representation is intransitive. One more example: the $n$-orbit is
combinatorially a $|F|\times n$ matrix in which are arranged $n$ elements
of a base set $V$.  This combinatorics leads to invariants of the group
that are not existing in group algebra. The study of these invariants
gives a hope to obtain a simple full invariant for a big class of groups.
For example it is of interest the next

\begin{hypothesis}\label{Full.Inv}
A primitive ld-group is defined by its order and degree.

\end{hypothesis}

or

\begin{problem}\label{Full.Inv}
Are there existing two primitive md-groups of a degree $n$ with the same
order but non-isomorphic $n$-orbits?

\end{problem}

For imprimitive groups there exist many such examples. We can also note
that a nmd-group is always imprimitive.

Now we consider one property that is related with a structure of
permutation groups. In \cite{Cameron} is described a hypothesis
(polycirculant conjecture) of M.~Klin and D. Maru$\breve{\rm
s}$i$\breve{\rm c}$ that the automorphism group of a transitive graph
contains a regular element (a permutation decomposed in cycles of the same
length). On group theory language this conjecture has the next
interpretation.

Let $F$ be a finite group and $A<F$ be a md-stabilizer. Let $F$ contains a
subgroup $P$ of a prime order $p$ that is intersected with no subgroup
from stabilizers of $F$ conjugated with $A$, then it follows that $P$ is
a regular subgroup of a representation $F(F/A)$. So the polycirculant
conjecture statements that if $F(F/A)$ is a graph automorphism group, then
$F$ contains the corresponding subgroup $P$.

The attempts to use such approach was not successful. The reason is that
the such interpretation of the problem is out of the inside structure of a
$n$-orbit, but the desired property is directly connected with symmetries
of $n$-orbits of corresponding permutation groups (s.\cite{Golubchik}).

It is known the existence of external and internal automorphisms of finite
groups. The similar property also exists in $n$-orbits. If $X_n$ is a
$n$-orbit of a permutation group $G$ of a degree $n$, then it can contain
isomorphic $k$-subsets that are connected with permutation of $G$ or with
permutation of symmetric group $S_n$ (the simplest example of the latter
is given by intransitive groups). This property of a $n$-orbit has only
partial interpretation in group theory, but it plays an important role in
construction of a $n$-orbit and hence in construction of the related
group.

The next is a combinatorial property of some normal subgroups.

\begin{lemma}\label{md-impr}
Let $G(V)$ be an imprimitive md-group and $Q$ be a partition of $V$ on
imprimitivity blocks (i.e. the action of $G$ on $V$ maintains $Q$), then
the maximal subgroup $Stab(\overrightarrow{Q})<G$ that maintains all
classes of $Q$ (i.e. arbitrarily ordered set $Q$) is a normal subgroup of
$G$.

\end{lemma}
\proof
Subgroups $\{Stab(U)<G:\, U\in Q\}$ form a class of conjugate subgroups
and $Stab(\overrightarrow{Q})=\cap_{(U\in Q)}Stab(U)$. For a md-group
$Stab(\overrightarrow{Q})$ is not trivial, because on definition an
imprimitive md-group has not isomorphic representation on $Q$.
$\Box$\bigskip

From here immediately follows that a md-representation of a simple group is
primitive. A simple proof of the Feit-Thompson theorem follows then from
the following statement.

\begin{theorem}\label{pr-even}
Let $G(V)$ be a primitive group of odd order, then it contains a normal
subgroup.

\end{theorem}

Now we shall prove this theorem.

\begin{lemma}\label{1-fixed}
Any permutation from $G$ fixes at most one element of $V$.

\end{lemma}
\proof
First, we can assume that no subgroup of $G$ that has a primitive
representation contradicts with lemma. Second, if a permutation $g$ of $G$
fixes $k$ element of $V$ (i.e. its decomposition in cycles contains
exactly $k$ cycles of the length $1$), then $k$ is an \emph{automorphic
number}, i.e.  $G$ contains a $k$-element \emph{suborbit} (or some subgroup
of $G$ contains a $k$-element orbit).  Moreover, if $g$ fixes $k$-tuple
$\alpha_k$, then a set $Co(\alpha_k)$ of coordinates of $\alpha_k$ is a
suborbit of $G$.

Let for $G$ the lemma be not correct and $k>1$ be a number of elements of
$V$ which are fixed by permutation $g$ of $G$. Let $I_k$ be fixed by $g$
$k$-subspace ($k$-tuple), $X_n$ be a $n$-orbit of $G$ and $X_k$ be a
projection of $X_n$ on $I_k$ ($X_k=\hat{p}(I_k)X_n=GI_k$), then $|X_k|>n$,
so there exist permutations of $G$ that fix only one element of $V$.

Let a stabilizer $Stab(v_1)$ of some element $v_1\in V$ has non-regular
transitive representation $A(U^1)$ on subset $U^1$ of $V$ (i.e.
$|A(U^1)|>|U^1|$), then a stabilizer of $A(U^1)$ fixes certain odd number
of elements from $U^1$.  So in order to fix odd number of elements of $V$
with permutations of $Stab(v_1)$ there must exist a subdirect product
$A(U^1)\circ\ldots\circ A(U^l)$ on an even number of non-intersected
subsets $U^1,\ldots,U^l$, where $\{v_1\}\cup U^1\cup\ldots\cup U^l=V$.
Then there exists only two possibilities. First, subsets $U^1,\ldots,U^l$
are $G$-isomorphic, then $G$ contains a subgroup whose projection on
$V\setminus \{v_1\}$ is imprimitive and hence $G$ contains an involution.
Contradiction. In second we consider a fixed $k$-tuple $\alpha_k$ that
contains all fixed elements from $U^1,\ldots,U^l$ and $v_1$. A stabilizer
$B=Stab(\alpha_k)$ has non-trivial normalizer $N=Stab(Co(\alpha_k))$. From
here follows that projection of $B$ on $V\setminus \{v_1\}$ has
non-trivial normalizer that is a projection of $N$ on $V\setminus
\{v_1\}$. It is only possible if for $T^1=Co(\alpha_k)\cap U^1$ a subset
$Z^1=U^1\setminus T^1$ is automorphic.  Since $|Z^1|$ is even then we have
the contradiction again.  $\Box$\bigskip

\begin{corollary}\label{|Stab(v)|<n}
A stabilizer $A<G$ of an element $v\in V$ has regular representation on a
subset of $V$ and hence $|A|<n$.

\end{corollary}

\begin{corollary}\label{tr.subgr}
$G$ contains a regular subgroup $H$.

\end{corollary}
\proof
Stabilizers of $G$ have trivial intersections, hence $G$ contains $n-1$
element which belong to no stabilizer. All this elements are regular,
because if one of those is not regular, then it generates a permutation
that fixes more than one element of $V$. $\Box$\bigskip

\begin{corollary}\label{ld}
$G$ is ld-group.

\end{corollary}
\proof
It follows from obtained above the structure of $G$. $\Box$\bigskip

So we have that $|H|>|A|$ and hence a representation of $G$ on classes of
$G/H$ is a homomorphism. It follows that $H$ is a normal subgroup of $G$.

\begin{corollary}\label{div.by.$4$}
Let $F$ be a simple group, then its order is divisible by $4$.

\end{corollary}
\proof
If $|F|=2m$ and $m$ is odd, then $F$ has evidently a normal subgroup.
$\Box$\bigskip


\section*{Conclusion}

A specificity of the $n$-orbit representation is a possibility to do a
group visible. In order to go to this visibility one makes a partition of
a matrix of a $n$-orbit on cells of $k$-orbits of subgroups of an
investigated permutation group and studies symmetry properties of cells
and the whole partition. It is not difficult to find that there exists
only three kinds of basic cells which allow to construct a $n$-orbit of
any group. Namely this approach to an investigation of permutation groups
gives the progress in some cases that are difficult for other existing
methods.

The using of $k$-orbits symmetry properties to the polynomial solution of
the graph isomorphism problem was begun by Author in 1984. The
generalization of $k$-orbits (regular $k$-sets) was used for describing of
the structure of strongly regular graphs and their generalization on
dimensions greater than two.

In 1997 Author understood the connection between the graph isomorphism
problem and the problem of a full invariant of a finite group and did
some attempts to obtain this full invariant by constructing of some
appropriate group representations. But again the best approach was
obtained with $k$-orbit representation.

Then this method was applied to the polycirculant conjecture and again with
success.

Thus all these themes wait for a possibility to be published.

Author will be thankful for any contact related with mentioned themes and
problems that are difficult in finite group theory and permutation group
theory.


\section*{Acknowledgements}
I would like to express many thanks to Prof.~P.~Cameron for his excellent
site on the internet, for two years of e-mail contacts, for presenting me
some examples of permutation groups interesting for analysis by methods
of $k$-orbit theory, and for pointing out mistakes in my first version
of this text.


\end{document}